\documentclass[reqno,times]{amsart}
\usepackage[english]{babel}  \usepackage{latexsym,amsfonts} \usepackage{amsmath, amssymb, amsthm} \usepackage{mathrsfs}
\usepackage{enumerate} \usepackage{lmodern} \usepackage{mathpazo}   \usepackage{euscript}  \usepackage[pdftex]{hyperref}

\usepackage[top=2cm,bottom=2cm,left=3cm,right=3cm]{geometry}
\numberwithin{equation}{section}  \makeatletter\@addtoreset{equation}{section}
\newtheorem {theorem}{Theorem}[section]        
  
   \newtheorem {corollary}[theorem]{Corollary}     \newtheorem {rem}[theorem]{Remark}
\newtheorem {proposition}[theorem]{Proposition}       
         
    \newcommand{\R}{\mathbb R}    	
\newcommand{\fin}{\hfill $\square$}

\begin{document}

\title{On a class of two-index real Hermite polynomials}
  \author{Naima A\"{i}t Jedda $\&$ Allal Ghanmi}
  \address{Department of Mathematics, P.O. Box 1014,  Faculty of Sciences, Mohammed V-Agdal University, Rabat, Morocco}
      \email{ag@fsr.ac.ma}  
\begin{abstract}
We introduce a class of doubly indexed real Hermite polynomials and we deal with their related properties
like the associated recurrence formulae, Runge's addition formula, generating function and Nielsen's identity.
\end{abstract}

\maketitle

\section{Introduction}
The Burchnall's operational formula (\cite{Burchnall41})
   \begin{align}\label{ROpfsum}
   \left(-\frac{d}{dx}+2x\right)^{m}(f) =
   m! \sum\limits_{k=0}^{m} \frac{(-1)^{k}}{k!} \frac{H_{m-k}(x)}{(m-k)!}\dfrac{d^{k}}{dx^{k}}(f),
   \end{align}
where  $H_{m}(x)$ denotes the usual Hermite polynomial (\cite{Hermite1864-1908,Rainville71})
   \begin{align}\label{Rodrigues}
   H_{m}(x)=(-1)^{m}e^{x^2}\dfrac{d^{m}}{dx^{m}}\left(e^{-x^2}\right),
  \end{align}
enjoy a number of remarkable properties. It is used by Burchnall \cite{Burchnall41} to give a direct
proof of Nielsen's identity (\cite{Nielsen18}) 
   \begin{align}\label{Nielsen}
   H_{m+n}(x) =m!n!\sum\limits_{k=0}^{min(m,n)}\frac{(-2)^k}{k!} \frac{H_{m-k}(x)}{(m-n)!}\frac{H_{n-k}(x)}{(n-k)!}.
   \end{align}
The special case of \eqref{ROpfsum} where $f=1$, i.e.,
   \begin{align}\label{Rel-Op}
   H_{m}(x)=\left(-\frac{d}{dx}+2x\right)^{m}\cdot (1).
   \end{align}
can be employed to recover in a easier way the generating function
   \begin{align}\label{gen-m}
   \sum\limits_{m=0}^{+\infty}H_{m}(x) \frac{t^{m}}{m!}=\exp (2xt-t^{2})
   \end{align}
as well as the Runge addition formula (\cite{Runge14,Kampe23})
  \begin{align} \label{addition-formula}
  H_{m}(x+y)=\left(\frac 12\right)^{m/2} m! \sum\limits_{k=0}^{n}\frac{H_{k}(\sqrt{2}x)}{k!} \frac{H_{m-k}(\sqrt{2}y)}{(m-k)!} .
  \end{align}

In this note, we have to consider the following class of doubly indexed real Hermite polynomials
  \begin{align}\label{DoubluIHP}
  H_{m,n}(x)=\left(-\frac{d}{dx}+2x\right)^{m}\cdot (x^n),
  \end{align}
and we derive some of their useful properties. More essentially, we discuss the associated
recurrence formulae, Runge's addition formula, generating function and Nielsen's identity.

 \section{  Doubly indexed real Hermite polynomials $H_{m,n}(x)$}

By taking $f(x)=x^n$ in \eqref{ROpfsum}, we obtain
   \begin{align}
   H_{m,n}(x)&:=\left(-\frac{d}{dx}+2x\right)^{m}(x^n) \label{hmn-def} \\
             & =m!n! \sum\limits_{k=0}^{min(m,n)} \frac{(-1)^k}{k!}
               \frac{x^{n-k}}{(n-k)!}\frac{H_{m-k}(x)}{(m-k)!}. \label{hmn}
    \end{align}
It follows that $H_{m,n}(x)$ is a polynomial of degree $m+n$, since
   $$
   Q(x):= H_{m,n}(x) - x^{n}H_{m}(x)
   $$
is a polynomial of degree $deg(Q)\leq n+m-2$. For the unity of the formulations, we shall define trivially
   $$
   H_{m,n}(x) = 0
   $$
whenever $m<0$ or $n<0$. We call them doubly indexed real Hermite polynomials. Note that $H_{m,0}(x)=H_{m}(x)$,
$H_{0,n}(x)=x^n$ and
   \begin{align} \label{hmn0}
    H_{m,n}(0) = \left\{ \begin{array}{lll} 0 & m < n \\
               (-1)^n  \frac{m! }{(m-n)!} H_{m-n}(0)   & m \geq n
             \end{array} \right. .
   \end{align}
A direct computation using  \eqref{hmn-def} gives rise to
   $$
   H_{1,n}(x) = -n x^{n-1} + 2x^{n+1}
   $$
 for every integer $n\geq 1$.  Note also that, since $ H_1(x)= 2x$, it follows
   \begin{align}\label{Hm1}
   H_{m+1}(x) = \left(-\frac{d}{dx}+2x\right)^{m}(H_1(x))
              = \left(-\frac{d}{dx}+2x\right)^{m}(2x)
              = 2H_{m,1}(x)  .
   \end{align}
The first few values of $H_{m,n}$ are given by

   \hspace*{.2cm}
 \begin{center}
 \begin{tabular}{||c||c|c|c||}
   \hline
   $H_{m,n}$ & $n=1$                   & $n=2$                   &  $n=3$          \\
   \hline\hline
   $m=1$     & $-1+2x^2 $              & $-2x +2x^3 $            &  $-3x^2 +2 x^4$ \\
   \hline
   $m=2$     & $-6x+4x^3 $             & $2 -10 x^2 + 4 x^4$     &  $6x -14 x^3 + 4 x^5$  \\
   \hline
   $m=3$     & $6 -24 x^2 + 8x^4 $     & $24x -36 x^3 + 8 x^5$   &  $-6 +54 x^2 - 48 x^4 + 8x^6$ \\
   \hline
 \end{tabular}
 \end{center}

 \hspace*{.2cm}

\noindent From \eqref{hmn}, one can deduce easily the symmetry formula
  \begin{align}
  H_{m,n}(-x) =(-1)^{n+m}H_{m,n}(x), \label{even-odd}
  \end{align}
so that the $H_{m,n}(x)$ is odd (rep. even) if and only if $n+m$ is odd (resp. even).
Furthermore, let mention that the Rodrigues formula for  $H_{m,n}(x)$ reads
   \begin{align} \label{RodriguesF}
   H_{m,n}(x)=(-1)^{m}e^{x^2} \frac{d^m}{dx^m} \left(x^{n}e^{-x^2}\right).
   \end{align}
Indeed, this is evidently proved using
   \begin{align}\label{ROpf}
   \left(-\frac{d}{dx}+2x\right)^{m}\cdot(f)=(-1)^{m}e^{x^2}\frac{d^m}{dx^m}\left(e^{-x^2}f\right)  .
   \end{align}
Therefore, these polynomials constitute a subclass of the generalized Hermite polynomials
   \begin{align} \label{GH-HP}
    H_m^\gamma(x,\alpha,p): = (-1)^m x^{-\alpha} e^{px^\gamma}
    \frac{d^m}{dx^m}\left(x^{\alpha} e^{-px^\gamma} \right).
    \end{align}
considered by Gould and Hopper in \cite{GouldHopper62}. In fact, we have $H_{m,n}(x)= x^n H_m^2(x,n,1)$.

    \begin{proposition} \label{prop-RF}
    The polynomials $H_{m,n}$; $m,n\geq 1$, satisfy the following recurrence formulae
    \begin{align}
    & H_{m,n}'(x)+ H_{m+1,n}(x) - 2xH_{m,n}(x) =0, \label{prop-RF1} \\
    & H_{m,n}(x) + nH_{m-1,n-1}(x)    - 2H_{m-1,n+1}(x)=0, \label{prop-RF2} \\
    & H_{m,n}(x) + m H_{m-1,n-1}(x) - x H_{m,n-1}(x)=0, \label{prop-RF3}\\
    & (m-n) H_{m-1,n-1}(x) + 2H_{m-1,n+1}(x) + x H_{m,n-1}(x)=0. \label{prop-RF4}
    \end{align}
    \end{proposition}

    \noindent {\it Proof.}
    The first one follows by writing the derivation operator as
    $$
    \frac{d}{dx}=-\left(-\frac{d}{dx}+2x\right) + 2x.
    $$
    Indeed, we get
   \begin{align*}
   \frac{d}{dx}\left(H_{m,n}(x)\right)
        &=-\left(-\frac{d}{dx}+2x\right)H_{m,n}(x) + 2xH_{m,n} (x)\\
        &=-H_{m+1,n}(x) + 2xH_{m,n}(x).
   \end{align*}
   For the second one, one writes $H_{m,n}(x)$ as
    \begin{align*} H_{m,n}(x)
                          &=\left(-\frac{d}{dx}+2x\right)^{m-1}\left(H_{1,n}(x)\right)\\
                          &=\left(-\frac{d}{dx}+2x\right)^{m-1}\left(-nx^{n-1}+2x^{n+1}\right)\\
                          &=-nH_{m-1,n-1}(x)+2H_{m-1,n+1}(x).
                          \end{align*}
To prove \eqref{prop-RF3}, we use \eqref{RodriguesF} combined with Leibnitz formula. Indeed,
  \begin{align*} H_{m,n}(x) &=(-1)^{m}e^{x^2} \frac{d^m}{dx^m} \left(x \cdot x^{n-1}e^{-x^2}\right) \\
                          &=(-1)^{m}e^{x^2}\left[ x \frac{d^m}{dx^m} \left(x^{n-1}e^{-x^2}\right) + m \frac{d^{m-1}}{dx^{m-1}} \left(x^{n-1}e^{-x^2}\right)\right] \\
                          &=   x H_{m,n-1}(x) - m H_{m-1,n-1}(x).
                          \end{align*}
  Finally,  \eqref{prop-RF4} follows from  \eqref{prop-RF2} and  \eqref{prop-RF3} by substraction.
    \fin

\begin{rem}
According to \eqref{Hm1}, the \eqref{prop-RF3}  (corresponding to $n=1$)
leads to the well known recurrence formula $H_{m+1}(x) = 2xH_m(x) - 2mH_{m-1}(x)$ for $H_m(x)$. Note also that
\eqref{prop-RF1} reduces further to $H_{m}'(x) + H_{m+1}(x) - 2xH_{m}(x) =0$ by taking $n=0$.
\end{rem}

  \begin{proposition} We have the following addition formula
      \begin{align} \label{add-forRealDHP} H_{m,n}(x+y)=m!n!\left(\frac 1{\sqrt{2}}\right)^{m+n} \sum\limits_{k=0}^{m} \sum\limits_{j=0}^{n}
     \frac{H_{k,j}(\sqrt{2}x)}{k!j!} \frac{H_{m-k,n-j}(\sqrt{2}y)}{(m-k)!(n-j)!} .
      \end{align}
 \end{proposition}

 \noindent {\it Proof.}   We have
  \begin{align*}
  H_{m,n}\left( x+y \right)  & =  \left(-\frac{d}{d(x+y)}+ 2(x+y)\right)^m.((x+y)^n)\\
            & = \left(-\frac 12 \left(\frac{\partial}{\partial x}+\frac{\partial}{\partial y}\right)+ 2(x+y)\right)^m.((x+y)^n)\\
            & = \left(\frac1{\sqrt{2}}\right)^{m} \left( A_x +A_y\right)^m .((x+y)^n)\\
             & = \left(\frac1{\sqrt{2}}\right)^{m} \sum_{j=0}^n     \binom{n}{j}     \left( A_x +A_y\right)^m .(x^jy^{n-j}),
  \end{align*}
 where $A_t$ stands for $ A_t = -{\partial}/{(\partial \sqrt 2t)}+ 2\sqrt 2 t.$
Thus, since $A_x$ and $A_y$ commute, we can make use of the binomial formula to get
                     \begin{align*}
  H_{m,n}\left( x+y \right)
            & = \left(\frac1{\sqrt{2}}\right)^{m} \sum\limits_{k=0}^{m} \sum_{j=0}^n      \binom{m}{k}         \binom{n}{j}     A_x^k.(x^j)
                 A_y^{m-k}.(y^{n-j}),
 \end{align*}
whence, we obtain the asserted result according to the fact that $$A_t^r(t^s) = {2}^{-s/2}H_{r,s}(\sqrt{2}t).$$
  \fin

  \begin{rem}
 We recover the Runge addition formula \eqref{addition-formula} for the classical real Hermite polynomials $H_m(x)=H_{m,0}(x)$ by taking $n=0$ in
 \eqref{add-forRealDHP}.
\end{rem}

 The following identities are immediate consequence of the previous proposition.
 \begin{corollary} The identity
\begin{align*}
H_{m,n}(t) =m!n!\left(\frac 1{\sqrt{2}}\right)^{m+n}  \sum\limits_{j=0}^{n} \sum\limits_{k=j}^{m}
\frac{(-1)^{j}}{j!(k-j)!} H_{k-j}(0) \frac{H_{m-k,n-j}(\sqrt{2}t)}{(m-k)!(n-j)!}
     \end{align*}
holds by taking $x=0$ and setting  $t=y$ in \eqref{add-forRealDHP}, keeping in mind \eqref{hmn0}. We get also
\begin{align*}  H_{m,n}(t)=m!n!\left(\frac 1{\sqrt{2}}\right)^{m+n} \sum\limits_{k=0}^{m} \sum\limits_{j=0}^{n}
     \frac{H_{k,j}( t/{\sqrt{2}})}{k!j!} \frac{H_{m-k,n-j}( t/{\sqrt{2}})}{(m-k)!(n-j)!} \end{align*}
     by setting $x=y=t/2$ in \eqref{add-forRealDHP}. While for $t=-\sqrt{2}x=\sqrt{2}y$, we obtain
\begin{align*}    \sum\limits_{k=0}^{m} \sum\limits_{j=0}^{n}
     (-1)^{k+j} \frac{H_{k,j}(t)}{k!j!} \frac{H_{m-k,n-j}(t)}{(m-k)!(n-j)!} = 0 \end{align*}
     whenever $m+n$ is odd or $m >n$.
\end{corollary}

Next, we state the following

   \begin{proposition}
   The generating function of $H_{m,n}$ is given by
    \begin{align} \label{gen-mn}
    \sum\limits_{m,n=0}^{+\infty}H_{m,n}(x)\frac{u^m}{m!}\frac{v^n}{n!}=\exp\left(-u^2+(2u+v)x-uv\right).
    \end{align}
  \end{proposition}

   \noindent {\it Proof.}  According to the definition of $H_{m,n}$, we can write
    \begin{align*}
    \sum\limits_{m,n=0}^{+\infty} H_{m,n}(x)\frac{u^m}{m!}\frac{v^n}{n!}
         & =    \left[\sum\limits_{m=0}^{+\infty} \frac{1}{m!} \left(-u\frac{d}{dx}+2ux\right)^m \right]\cdot
                \left(\sum\limits_{n=0}^{+\infty}\frac{v^n}{n!}x^n\right)\\
         & =    \exp\left( -u\frac{d}{dx}+2ux\right) \left(e^{vx}\right).
  \end{align*}
  Making use of the Weyl identity which reads for the operators $A = 2 xId $ et $B  = -{d}/{dx}$ as
                \begin{align*} \label{BCH}
              \exp(uA+uB)= \exp(uA) \exp(uB)\exp\left(-u^2Id\right); \quad u\in \R,
              \end{align*}
   we get
      \begin{align*}
    \sum\limits_{m,n=0}^{+\infty} H_{m,n}(x)\frac{u^m}{m!}\frac{v^n}{n!}
       & =    e^{2ux -u^2} \exp\left( -u\frac{d}{dx}\right) \left(e^{vx}\right).
  \end{align*}
  Therefore, the desired result follows since
     $$
     \exp\left( -u\frac{d}{dx}\right) \left(e^{vx}\right)
     = \sum\limits_{k=0}^\infty  \frac{(-u)^k}{k!}\left(\frac{d}{dx}\right)^k (e^{vx})
     =    e^{-uv} e^{vx}.
     $$
   \fin

 \begin{rem}
The special case of $v=0$ (in \eqref{gen-mn}) infers the generating function \eqref{gen-m} of the standard real Hermite polynomials $H_m$.
Furthermore, for $y=u=-v$, we get
\begin{equation}
   e^{xy}=  \sum\limits_{m,n=0}^{+\infty}(-1)^n H_{m,n}(x)\frac{y^{m+n}}{m! n!}. \label{gen-mnuv}
     \end{equation}
       \end{rem}

  \begin{proposition}
   We have the recurrence formula
    \begin{equation}
      H_{m,n}'(x) = 2m H_{m-1,n}(x) + n H_{m,n-1}(x) .\label{gen-mm}
     \end{equation}
  \end{proposition}

 \noindent {\it Proof.}
  Differentiating the both sides of \eqref{gen-mn} and making appropriate changes of indices yield \eqref{gen-mm}.
 \fin

\begin{corollary}\label{corDerHmn} We have
   \begin{equation}\label{derHmn}
   \frac{d^\nu}{dx^\nu}(H_{r,n}(x)) = r!n! \sum_{j=0}^\nu \alpha_{j,\nu}  \frac{H_{r-\nu+j,n-j}(x)}{ (r-\nu+j)!(n-j)! }  ,
   \end{equation}
     where
     $$ \alpha_{j,\nu} =  \left\{
                                \begin{array}{lll}
                                  2^\nu  & \quad \mbox{for} \quad  j=0 \\
                                   2 \alpha_{j,\nu-1} + \alpha_{j-1,\nu-1}   & \quad \mbox{for} \quad  1 \leq j < \nu \\
                                  1    & \quad \mbox{for} \quad j=\nu \\
                                \end{array}
                              \right.
     .$$
\end{corollary}

 \noindent {\it Proof.}
 This can be handled by mathematical induction using \eqref{gen-mm}.
  \fin

 \begin{rem}
 The $\alpha_{j,\nu}$ are even positive numbers and their first values are
 $$
   \begin{array}{c||cccccc}
   \alpha_{j,\nu} &  j=0     &    j=1    &   j=2     &    j=3   &    j=4    & j=5 \\
\hline\hline
   \nu=0          &        1 &           &           &          &           &     \\
   \nu=1          & \fbox{2} &       1   &           &          &           &     \\
   \nu=2          &      2^2 & \fbox{4}  &       1   &          &           &     \\
   \nu=3          &      2^3 &      12   & \fbox{6}  &      1   &           &     \\
   \nu=4          &      2^4 &      32   &      24   & \fbox{8} &     1     &     \\
   \nu=5          &      2^5 &      80   &      80   &      40  & \fbox{10} & 1   \\
   \end{array}
.
 $$
 \end{rem}

 We conclude this note by giving a formula for the two-index Hermite polynomial $H_{m,n}(x)$
expressing it as a weighted sum of a product of the same polynomials. Namely, we state the following

   \begin{proposition}
   Keep notation as above. Then the Nielsen identity for $H_{m,n}$; $n\geq 1$,  reads
       \begin{align*}
       H_{m+r,n}(x) &=    m! r!n n! \sum_{k,\nu,j=0}^{m,k,\nu} \alpha_{j,\nu} \frac{\Gamma(n+k-\nu)}{(k-\nu)! \nu!}
    \frac{(-x)^{\nu}}{x^{n+k}} \frac{H_{m-k,n}(x)}{(m-k)!n!}    \frac{H_{r-\nu+j,n-j}(x)}{ (r-\nu+j)!(n-j)! }.
   \end{align*}
  \end{proposition}

 \noindent {\it Proof.}  Recall first that $H_m^\gamma(x,\alpha,p)$, the polynomials given through
\eqref{GH-HP}, can be rewritten in the following equivalent form (\cite{GouldHopper62})
  $$
  H_{m}^\gamma(x,\alpha,p) := \left( -\frac{d}{dx} +p\gamma x^{\gamma-1} - \frac \alpha x\right)^m(1).
  $$
Now, since for the special values $p=1$, $\gamma= 2$ and $\alpha=n$, we have
\begin{align*}
H_{m+r,n}(x) &=   x^n H_{m+r} ^2(x,n,1) \\
&=   x^n \left( -\frac{d}{dx} + 2 x  - \frac n x\right)^{m}\left(H_{r} ^2(x,n,1)\right)\\
&=   x^n \left( -\frac{d}{dx} + 2 x  - \frac n x\right)^{m}\left(x^{-n} H_{r,n}(x) \right),
\end{align*}
we can make use of the Burchnall's formula extension proved by Gould and Hopper \cite{GouldHopper62}, to wit
 $$
 \left(
 -\frac{d}{dx} +p\gamma x^{\gamma-1} - \frac \alpha x\right)^m(f)
 = m!\sum_{k=0}^m \frac{(-1)^{k}}{k!} \frac{H_{m-k}^\gamma(x,\alpha,p)}{(m-k)!} \frac{d^k}{dx^k}(f)
. $$
Thus, for $f= x^{-n}H_{r,n}$, we obtain
\begin{align}
 H_{m+r,n}(x)
 &= m! \sum_{k=0}^m \frac{(-1)^{k}}{k!} \frac{H_{m-k,n}(x)}{(m-k)!} \frac{d^k}{dx^k}(x^{-n}H_{r,n}(x))
 .\label{nielsen-n} \end{align}
 Therefore,  by applying the Leibnitz formula and appealing the result of Corollary \ref{corDerHmn}, we get
 \begin{align*}
 H_{m+r,n}(x)
  &= m! \sum_{k=0}^m \frac{(-1)^{k}}{k!} \frac{H_{m-k,n}(x)}{(m-k)!}
  \sum_{\nu=0}^k \binom{k}{\nu} \frac{d^{k-\nu}}{dx^{k-\nu}}(x^{-n}) \frac{d^\nu}{dx^\nu}( H_{r,n}(x))\\
  &{=}
   m! r!n n! \sum_{k,\nu,j=0}^{m,k,\nu}  
    \alpha_{j,\nu} \frac{\Gamma(n+k-\nu)}{(k-\nu)! \nu!}
    \frac{(-x)^{\nu}}{x^{n+k}} \frac{H_{m-k,n}(x)}{(m-k)!n!}    \frac{H_{r-\nu+j,n-j}(x)}{ (r-\nu+j)!(n-j)! }
 \end{align*}
 for every integer $n\geq 1$. Note that for $n=0$, \eqref{nielsen-n} reads simply
 $$
 H_{m+r}(x)
 = m! \sum_{k=0}^m \frac{(-1)^{k}}{k!} \frac{H_{m-k}(x)}{(m-k)!} \frac{d^k}{dx^k}(H_{r}(x)).
 $$
In this case,  we recover the usual Nielsen formula \eqref{Nielsen} for the real Hermite polynomials $H_m$.
 \fin

\end{document}